# ON PATHWISE UNIQUENESS FOR REFLECTING BROWNIAN MOTION IN $C^{1+\gamma}$ DOMAINS[1]

By Richard F. Bass and Krzysztof Burdzy

*University of Connecticut and University of Washington*

Pathwise uniqueness holds for the Skorokhod stochastic differential equation in $C^{1+\gamma}$ domains in $\mathbb{R}^d$ for $\gamma > 1/2$ and $d \geq 3$.

**1. Introduction.** We will prove pathwise uniqueness for the Skorokhod equation in every $C^{1+\gamma}$ domain $D \subset \mathbb{R}^d$, with $\gamma > 1/2$ and $d \geq 3$. We begin by giving rigorous definitions of these terms and state our main result.

We will write $x = (x_1, \ldots, x_{d-1}, x_d) = (\hat{x}, x_d)$ for $x \in \mathbb{R}^d$. We will say that a function $\Phi : \mathbb{R}^{d-1} \to \mathbb{R}$ is $C^{1+\gamma}$ if $\Phi$ is bounded and for some constant $c_1$ and all $\hat{x}$ and $\hat{y}$,

$$|\nabla\Phi(\hat{x}) - \nabla\Phi(\hat{y})| \leq c_1 |\hat{x} - \hat{y}|^\gamma, \qquad \hat{x}, \hat{y} \in \mathbb{R}^{d-1}.$$

A $C^{1+\gamma}$ domain is one which can be represented locally as the region above the graph of a $C^{1+\gamma}$ function.

Reflecting Brownian motion in a Euclidean domain with Lipschitz boundary can be represented as a solution to the following Skorokhod stochastic differential equation:

$$(1.1) \qquad X_t = X_0 + W_t + \int_0^t \mathbf{n}(X_s)\, dL_s,$$

where $\mathbf{n}(x)$ is the inward pointing unit normal vector at $x \in \partial D$, $W = \{W_t : t \geq 0\}$ is a $d$-dimensional Brownian motion with respect to a filtered probability space $(\Omega, \mathcal{F}, \{\mathcal{F}_t\}, \mathbb{P})$, $X_0 \in \overline{D}$, and we require $X_t \in \overline{D}$ for all $t$. Moreover, $L_t$ is a nondecreasing continuous process that increases only when $X_t \in \partial D$, that is, $\int_0^\infty \mathbf{1}_{\{X_t \in D\}}\, dL_t = 0$, a.s.

Received June 2007; revised October 2007.
[1]Supported in part by NSF Grants DMS-06-01783 and DMS-06-00206.

THEOREM 1.1. *If $X$ and $X'$ are two solutions to the Skorokhod equation with respect to the same Brownian motion (but with possibly different filtrations $\{\mathcal{F}_t\}$ and $\{\mathcal{F}'_t\}$, resp.) in a $C^{1+\gamma}$ domain $D \subset \mathbb{R}^d$, with $\gamma > 1/2$ and $d \geq 3$, then*

$$\mathbb{P}(X_t = X'_t \text{ for all } t \geq 0) = 1.$$

Weak existence and weak uniqueness for (1.1) are known; see, for example, Theorem 4.2 in [6]. This and the pathwise uniqueness proved in Theorem 1.1 imply strong existence by a standard argument (see [10], Theorem IX.1.7).

Pathwise uniqueness for (1.1) has been proved for all Lipschitz domains $D \subset \mathbb{R}^2$ with Lipschitz constant 1 in [5]. See [4] for a shorter proof. In particular, pathwise uniqueness holds in $C^{1+\gamma}$ domains in $\mathbb{R}^2$ for every $\gamma > 0$. See [3] for an introduction to the theory of reflecting Brownian motion and [5] for the history of the problem considered in this paper. It was asserted in [6] that pathwise uniqueness holds in $C^{1+\gamma}$ domains for all $\gamma > 0$. The proof given in that paper is incorrect; see Remark 3.8 for further details.

We have a heuristic argument showing that pathwise uniqueness fails in some $C^{1+\gamma}$ domains for some $\gamma > 0$ and some $d \geq 3$. At this time, we are not able to supply all the details needed to turn this claim into a rigorous proof. However, we present a counterexample showing that $\gamma = 1/2$ is the critical value for Proposition 3.4, a result on pathwise uniqueness for reflected SDEs in a half-space that are closely related to the reflecting Brownian motion in $C^{1+\gamma}$ domains. The counterexample gives credence to our belief that Theorem 1.1 is sharp, that is, it does not hold for $\gamma < 1/2$.

The proof of Theorem 1.1 is based on an analysis of the distance between two solutions to (1.1). We use Itô's formula to express the distance as the sum of a martingale and a process with finite variation. The size of the oscillations for each of these processes is bounded above using estimates for the Green function of the solution to (1.1) and bounds for partial derivatives of a mapping of $D$ to a half-space. Some of the crucial estimates are taken from the literature on PDEs in $C^{1+\gamma}$ domains. We conclude that the distance between the two solutions has to be 0, hence pathwise uniqueness holds.

The paper has three more sections. The next section gives some estimates for solutions to partial differential equations. The proof of the main theorem is given in the third section. A counterexample showing that Proposition 3.4 is sharp is given in the last section.

**2. Estimates for solutions to elliptic PDEs.** First, we introduce some notation. By standard techniques, it is enough to prove uniqueness in the region above the graph of a $C^{1+\gamma}$ function $\Phi$. Suppose $\Phi : \mathbb{R}^{d-1} \to \mathbb{R}$, and $D = \{(\hat{x}, x_d) : x_d > \Phi(\hat{x})\}$.

Let
$$Q(x,r) = \{y : |y_i - x_i| \leq r, i = 1, \ldots, d\},$$
$$B(x,r) = \{y \in \mathbb{R}^d : |y - x| < r\},$$
$$\hat{B}(\hat{x}, r) = \{\hat{y} \in \mathbb{R}^{d-1} : |\hat{y} - \hat{x}| < r\}.$$

Let $\Gamma : \mathbb{R}^d \to \mathbb{R}^d$ be defined by
$$\Gamma(x) = (\hat{x}, x_d - \Phi(\hat{x})).$$

We use $U = \{x_d > 0\}$ for the upper half-space.

Let
$$|\Phi|_{C^{1+\gamma}} = \sup_{\hat{x} \neq \hat{y}} \frac{|\nabla \Phi(\hat{x}) - \nabla \Phi(\hat{y})|}{|\hat{x} - \hat{y}|^\gamma},$$
the usual seminorm for the space $C^{1+\gamma}$.

We will use the well-known fact that if $h$ is harmonic in a ball $B(x, r)$ and bounded by $M$ in absolute value in $B(x, r)$, then

(2.1)
$$|\nabla h(x)| \leq \frac{c_1 M}{r};$$

see [1], Corollary II.1.4, for example. By induction, the $j$th-order partial derivatives of $h$ at $x$ are bounded in absolute value by $c_j M / r^j$.

PROPOSITION 2.1. *For any $c_0 < \infty$ there exist $c_1, c_2 < \infty$ such that the following is true. Suppose that $\Phi(\hat{0}) = 0$, $\nabla \Phi(\hat{0}) = 0$, $\sup_{\hat{x}} |\Phi(\hat{x})| \leq 1/4$, $|\Phi|_{C^{1+\gamma}} \leq c_0$, and that $\Phi(\hat{x}) = 0$ if $|\hat{x}| \geq 3/4$. Suppose $h$ is harmonic and bounded in absolute value by 1 in $D \cap Q(0, 1)$ and $h$ is $C^2$ on each face of $\partial(Q(0,1) \cap D) \setminus \partial D$. Suppose either (i) $h = 0$ on $Q(0,1) \cap \partial D$, or (ii) $\partial h / \partial \mathbf{n} = 0$ a.e. on $Q(0,1) \cap \partial D$. Then:*

(a) *we have*
$$|\nabla h(x) - \nabla h(y)| \leq c_1 |x - y|^\gamma, \qquad x, y \in B(0, 5/8) \cap \overline{D},$$
*and $|\nabla h(x)| \leq c_1$ for $x \in B(0, 5/8) \cap \overline{D}$;*

(b) *for any $i, j \leq d$, we have*
$$\left| \frac{\partial^2 h}{\partial x_i \partial x_j}(x) \right| \leq c_2 \operatorname{dist}(x, \partial D)^{\gamma - 1}, \qquad x \in B(0, 1/2) \cap D.$$

PROOF. (a) Case (i) follows immediately from Theorem A of [7], if we take $A(x, y, p) = p$ and $B(x, y, p) = 0$ in their theorem.

Let us consider case (ii). First we claim that $h$ is $C^\infty$ in $A_1 \cup A_2$, where
$$A_1 = (\hat{B}(\hat{0}, 15/16) \setminus \hat{B}(\hat{0}, 13/16)) \times (0, 7/8),$$
$$A_2 = \hat{B}(\hat{0}, 15/16) \times (1/2, 7/8).$$

Moreover, we claim that the absolute values of the $k$th partial derivatives of $h$ are bounded in $A_1 \cup A_2$ by constants depending only on $k$ and $c_0$ but otherwise independent of $\Phi$ and $h$.

The result for $A_2$ follows from (2.1) since $h$ is bounded in absolute value by 1 and $A_2$ is a positive distance from the boundary of $D$. We turn to $A_1$. Consider any vertical line segment in $A_1$ of the form $L = \hat{x}_0 \times (0, 7/8)$. Let $T = \hat{B}(\hat{x}_0, 1/32) \times (-15/16, 15/16)$ be a tube about this line segment and define $f(x)$ on $\partial T$ to be equal to $h(\hat{x}, |x_d|)$. If $Z$ is a standard $d$-dimensional Brownian motion and $\tau(T)$ is the first exit of this process from $T$, then the function

$$(2.2) \qquad \widetilde{h}(x) = \mathbb{E}^x f(Z_{\tau(T)}) + \mathbb{E}^{(\hat{x}, -x_d)} f(Z_{\tau(T)})$$

is easily seen to have the same boundary values as $h$ on $\overline{D} \cap \partial T$. The function $\widetilde{h}$ is harmonic because each term in (2.2) is separately harmonic in $T$. We have $\partial \widetilde{h}/\partial \mathbf{n} = 0$ on $T \cap \partial D$ by symmetry. Hence $\widetilde{h} = h$ in $T \cap U$. This implies that $h$ is $C^\infty$ on $A_1$. The estimate (2.1) and the remark following it can be applied to $\widetilde{h}$ inside $B(x, 1/32)$ for $x \in L$. We see that the absolute values of $k$th partial derivatives of $\widetilde{h}$ are bounded on $L$ by constants depending only on $k$. The same holds for $h$, because $h = \widetilde{h}$ on $L$. This finishes the proof of our claim about the behavior of $h$ and its derivatives on $A_1 \cup A_2$.

Since $h$ and $\partial h/\partial \mathbf{n}$ are bounded in $D_1 = (\hat{B}(\hat{0}, 7/8) \times (-7/8, 7/8)) \cap D$, by Green's first identity,

$$\int_{D_1} |\nabla h(x)|^2 \, dx = \int_{\partial D_1} h(x) \frac{\partial h}{\partial \mathbf{n}}(x) \rho(dx) < \infty,$$

where $\rho$ is surface measure on $\partial D_1$. Hence the function $h$ is in $W^{1,2}(D_1)$, the Sobolev space of functions whose gradient is in $L^2$. Choose a $C^\infty$ function $g : D \to \mathbb{R}$ of compact support that equals $h$ in $\overline{D_1} \cap D$ and has zero normal derivative on $\partial U \setminus D_1$. Such a $g$ is easy to find because $\partial D$ is flat outside of $D_1$.

Let $f = \Delta g$. We now apply Theorem 5.1 of [9] in the domain $\{x : x_d > \Phi(\hat{x})\}$. We set $\zeta(s) = s^\gamma$, so that $I(\zeta)(t) = \int_0^t \zeta(s) \frac{ds}{s} = c_3 t^\gamma$. We set $A(x, z, p) = p$ so that $\operatorname{div} A(x, u, \nabla u) = \operatorname{div} \nabla u = \Delta u$. We set $B(x, z, p) = -f(x)$, and therefore $u = g$ solves $\operatorname{div} A(x, u, \nabla u) + B(x, u, \nabla u) = 0$ in the region above the graph of $\Phi$. We set $\psi(x, u) \equiv 0$, so the boundary condition in the theorem in [9] becomes $\partial g/\partial \mathbf{n} = 0$. Set $\alpha = \gamma$. With these choices, verifying assumptions (5.1)–(5.6) of Theorem 5.1 in [9] is routine. According to a remark in the second to last paragraph on page 98 of [9], one can take $\beta = 1$ in Theorem 5.1 of that paper. The conclusion of part (a) of our lemma in case (ii) follows now from (5.7) in [9].

To prove (b), fix $i$ and let $g(x) = \frac{\partial h}{\partial x_i}(x)$. Then $g$ is harmonic in $D$ and is $C^\gamma$ in $B(0, 1/2) \cap D$. Fix $x \in B(0, 1/2) \cap D$ and let $r = \operatorname{dist}(x, \partial D)/4$. Let $\varphi$ be a nonnegative $C^2$ function supported in $B(0,1)$ with $\int \varphi = 1$, and let $\varphi_r(y) = r^{-d}\varphi(y/r)$. Let $g_r = g * \varphi_r$. Since

$$|g(y) - g_r(y)| = \left|\int [g(y) - g(y-z)]\varphi_r(z)\,dz\right|$$

$$\leq c_3 \int |z|^\gamma \varphi_r(z)\,dz \leq c_4 r^\gamma,$$

then by (2.1)

(2.3) $\quad |\nabla(g - g_r)(x)| \leq \dfrac{c_5}{r} \sup_{y \in \partial B(x,r)} |g(y) - g_r(y)| \leq c_6 r^{\gamma-1}.$

Note that $|\nabla \varphi_r(z)| \leq c_7 r^{-d-1}$. Using dominated convergence to justify interchanging differentiation and integration and the fact that $\int \varphi_r(z)\,dz = 1$, we have $\int \nabla \varphi_r(z)\,dz = 0$. Therefore

$$|\nabla g_r(x)| = \left|\int [g(x-z) - g(x)]\nabla \varphi_r(z)\,dz\right|$$

(2.4) $\quad\quad\quad \leq c_8 \int |z|^\gamma |\nabla \varphi_r(z)|\,dz$

$$\leq c_9 r^{\gamma-1}.$$

Combining (2.3) and (2.4), we obtain

$$|\nabla g(x)| \leq (c_6 + c_9) r^{\gamma-1},$$

which proves (b). $\square$

### 3. Pathwise uniqueness.

LEMMA 3.1. *Suppose $0 \in \partial D$ and $\Phi(\hat{0}) = 0$. There exist a continuous map $H = (h_1, \ldots, h_d) : D \to U$ with $H(0) = 0$, $\varepsilon > 0$ and $c_1, c_2 > 0$ such that:*

(a) *$H$ is one-to-one on $B(0,\varepsilon) \cap D$,*
(b) *$H$ maps $\partial D \cap B(0,\varepsilon)$ to $\partial U$,*
(c) *the functions $h_k$, $k = 1, \ldots, d$, are harmonic on $B(0,\varepsilon) \cap D$,*
(d) *$h_d(x) = 0$ for $x \in \partial D \cap B(0,\varepsilon)$,*
(e) *$\nabla h_k \cdot \mathbf{n} = 0$ on $\partial D \cap B(0,\varepsilon)$ for $k = 1, \ldots, d-1$,*
(f) *$\nabla h_d(x) = c(x)\mathbf{n}(x)$ for $x \in \partial D \cap B(0,\varepsilon)$,*
(g) *the Jacobian of $H$ is bounded above by 2 and below by $1/2$ on $B(0,\varepsilon) \cap D$,*
(h) *$\partial h_d / \partial x_d \geq c_1$ and*
(i) *for $i, j, k \leq d$,*

(3.1) $\quad \left|\dfrac{\partial^2 h_k}{\partial x_i \partial x_j}(x)\right| \leq c_2 \operatorname{dist}(x, \partial D)^{\gamma-1}, \quad x \in B(0,\varepsilon) \cap D.$

PROOF. It is routine to construct the following functions $\widetilde{\Phi}_\varepsilon$ and the corresponding domains $\widetilde{D}_\varepsilon$. For each $\varepsilon > 0$ let $\widetilde{\Phi}_\varepsilon$ agree with $\Phi$ on $\hat{B}(0,\varepsilon)$, $\widetilde{\Phi}_\varepsilon(\hat{x}) = 0$ if $|\hat{x}| \geq 3\varepsilon$, and $|\widetilde{\Phi}_\varepsilon|_{C^{1+\gamma}} \leq 2|\Phi|_{C^{1+\gamma}}$. Note that for some $\delta > 0$ and all $\varepsilon \in (0, \delta)$, $\|\widetilde{\Phi}_\varepsilon\|_\infty < \varepsilon$. Let

$$\widetilde{D}_\varepsilon = \{x : |x_i| < 8\varepsilon \text{ for } i = 1, \ldots, d-1 \text{ and } \widetilde{\Phi}_\varepsilon(\hat{x}) < x_d < 8\varepsilon\}.$$

Let $D_\varepsilon$ be the dilation of $\widetilde{D}_\varepsilon$ by the factor $(8\varepsilon)^{-1}$, that is, $D_\varepsilon = \{x : x/(8\varepsilon) \in \widetilde{D}_\varepsilon\}$. Let $\partial_L D_\varepsilon = \{x \in \partial D_\varepsilon : 8\varepsilon x_d = \widetilde{\Phi}_\varepsilon(8\varepsilon \hat{x})\}$.

For each $\varepsilon$, we define $h_i^\varepsilon$, $i \leq d$, on $D_\varepsilon$, as follows. For $i < d$, let $h_i^\varepsilon$ be the harmonic function whose boundary values are $x_i$ on $\partial D_\varepsilon \setminus \partial_L D_\varepsilon$, and whose normal derivative on $\partial_L D_\varepsilon$ is 0. Define $h_d^\varepsilon$ to be the harmonic function whose boundary values are $x_d$ on $\partial D_\varepsilon \setminus \partial_L D_\varepsilon$ and 0 on $\partial_L D_\varepsilon$.

By Proposition 2.1, we see that the $h_i^\varepsilon$ and their gradients are equicontinuous in $B(0, 1/2) \cap \overline{D}_\varepsilon$. Taking a sequence $\varepsilon_j \to 0$ along which the $h_i^{\varepsilon_j}$ and their gradients converge, it is clear that in the limit, we obtain for $i < d$ a function that is harmonic in $Q(0,1) \cap U$, has boundary values $x_i$ on $\partial Q(0,1) \cap U$, and has zero normal derivative on $Q(0,1) \cap \partial U$. Therefore the limit must be the function $x_i$. For $i = d$, a similar argument shows that the limit is the function $x_d$. The gradients also converge, so the limit of the gradients of $h_i^{\varepsilon_j}$ must be $e_i$, the unit vector in the $x_i$ direction.

Therefore, if we take $j$ sufficiently large, and set $\varepsilon = \varepsilon_j$, then the function

$$\overline{H}(x) = (h_1^{\varepsilon_j}(x/(8\varepsilon_j)), \ldots, h_d^{\varepsilon_j}(x/(8\varepsilon_j))), \qquad x \in B(0, \varepsilon_j) \cap D,$$

will have all the properties stated in the lemma on $B(0, \varepsilon_j) \cap D$. It is routine to find an extension $H$ of $\overline{H}$ to $D$ satisfying the desired properties on all of $D$. □

Let $E \subset \mathbb{R}^d$ be a domain and $Z$ a continuous process taking values in $\overline{E}$. Let us say $L_t^Z$ is a local time for $Z$ on $\partial E$ if $L_t^Z$ is a nondecreasing function that increases only when $Z \in \partial E$, that is, $\int_0^\infty \mathbf{1}_{\{Z_t \in E\}} dL_t^Z = 0$ a.s.

LEMMA 3.2. *Suppose that $H$ and $\varepsilon > 0$ are as in Lemma 3.1, $x_0 \in D \cap B(0, \varepsilon/2)$, and $X$ is a solution to*

$$X_t = x_0 + W_t + \int_0^t \mathbf{n}(X_s) \, dL_s^X,$$

*where $W = (W^1, \ldots, W^d)$ is a $d$-dimensional Brownian motion and $L_t^X$ is a local time of $X$ on $\partial D$. Assume that these processes are defined on a filtered probability space $(\Omega, \mathcal{F}, \{\mathcal{F}_t\}, \mathbb{P})$, and $W, X$ and $L^X$ are adapted to $\{\mathcal{F}_t\}$. Let $Y = (Y^1, \ldots, Y^d) = H(X)$ and let*

$$\tau = \inf\{t > 0 : X_t \notin D \cap B(0, 3\varepsilon/4)\}.$$

*Then for $t < \tau$ and $i = 1, \ldots, d-1$,*

$$Y_t^i - Y_0^i = \int_0^t \sum_{k=1}^d \frac{\partial h_i}{\partial x_k}(H^{-1}(Y_s))\, dW_s^k$$

*and*

$$Y_t^d - Y_0^d = \int_0^t \sum_{k=1}^d \frac{\partial h_d}{\partial x_k}(H^{-1}(Y_s))\, dW_s^k + \int_0^t \mathbf{n}(H^{-1}(Y_s))\, dL_s^Y,$$

*where $L_t^Y$ is a local time of $Y$ on $\partial U$ that is adapted to $\{\mathcal{F}_t\}$.*

PROOF. The lemma follows from Itô's formula, but we cannot apply Itô's formula directly because the function $H$ is not necessarily of class $C^2$ on $\overline{D}$. For $\delta > 0$, let

$$H_\delta((x_1, \ldots, x_d)) = H((x_1, \ldots, x_{d-1}, x_d + \delta)),$$

and $Y^\delta = (Y^{\delta,1}, \ldots, Y^{\delta,d}) = H_\delta(X)$. Note that for $\delta \in (0, \varepsilon/16)$, the components of $H_\delta$ are harmonic and $C^2$ on $\overline{D \cap B(0, 7\varepsilon/8)}$. Let $e_d = (0, \ldots, 0, 1)$ and $\mathbf{n}(x) = (\mathbf{n}_1(x), \ldots, \mathbf{n}_d(x))$, the inward pointing normal vector on $\partial D$. By Itô's formula, for $t < \tau$ and $i = 1, \ldots, d$,

(3.2)
$$Y_t^{\delta,i} - Y_0^{\delta,i} = \int_0^t \sum_{k=1}^d \frac{\partial h_i}{\partial x_k}(X_s + \delta e_d)\, dW_s^k$$
$$+ \int_0^t \sum_{k=1}^d \frac{\partial h_i}{\partial x_k}(X_s + \delta e_d)\mathbf{n}_k(X_s)\, dL_s^X.$$

By Proposition 2.1, the first partial derivatives of $h_i$'s are Hölder continuous with exponent $\gamma > 1/2$. This and the fact that $X$ is a continuous process imply that for any $i$ and $k$, the processes $\frac{\partial h_i}{\partial x_k}(X_s + \delta e_d)$ converge uniformly to $\frac{\partial h_i}{\partial x_k}(X_s)$ on the time interval $[0, \tau]$, a.s., when $\delta \to 0$. It follows that, a.s., when $\delta \to 0$,

(3.3) $$\int_0^t \sum_{k=1}^d \frac{\partial h_i}{\partial x_k}(X_s + \delta e_d)\, dW_s^k \to \int_0^t \sum_{k=1}^d \frac{\partial h_i}{\partial x_k}(X_s)\, dW_s^k.$$

By convention, let $\mathbf{n}(x)$ be the zero vector for $x \notin \partial D$. Recall from Lemma 3.1 that $\nabla h_k \cdot \mathbf{n} = 0$ on $\partial D \cap B(0, \varepsilon)$ for $k = 1, \ldots, d-1$ and $\nabla h_d(x) = c(x)\mathbf{n}(x)$ for $x \in \partial D \cap B(0, \varepsilon)$. We use the Hölder continuity of the $\nabla h_i$'s and the continuity of $X$ to conclude that for $i = 1, \ldots, d-1$ and all $k$,

(3.4) $$\frac{\partial h_i}{\partial x_k}(X_s + \delta e_d)\mathbf{n}_k(X_s) \to 0,$$

uniformly on $[0, \tau]$, a.s., when $\delta \to 0$. This also holds for $i = d$ and $k \neq d$. For $i = k = d$, we obtain

$$\frac{\partial h_i}{\partial x_k}(X_s + \delta e_d)\mathbf{n}_k(X_s) \to c(X_s)\mathbf{n}_k(X_s), \tag{3.5}$$

where $c(x)$ is a bounded function.

Since $H$ and $X$ are continuous, $Y_t^{\delta,i} - Y_0^{\delta,i} \to Y_t^i - Y_0^i$. We combine this observation with (3.2)–(3.5) to see that for $t < \tau$ and $i = 1, \ldots, d-1$,

$$Y_t^i - Y_0^i = \int_0^t \sum_{k=1}^d \frac{\partial h_i}{\partial x_k}(X_s)\, dW_s^k \tag{3.6}$$

and

$$Y_t^d - Y_0^d = \int_0^t \sum_{k=1}^d \frac{\partial h_d}{\partial x_k}(X_s)\, dW_s^k + \int_0^t \mathbf{n}(X_s)c(X_s)\, dL_s^X. \tag{3.7}$$

Note that $L_t^Y \stackrel{\mathrm{df}}{=} \int_0^t c(X_s)\, dL_s^X$ satisfies the definition of local time of $Y$ on $\partial U$, that is, it is a nondecreasing continuous process such that $\int_0^\tau \mathbf{1}_{\{Y_s \in D\}}\, dL_s^Y = 0$, a.s., and is adapted to $\{\mathcal{F}_t\}$. Since $H$ is one-to-one, we can replace $X_s$ with $H^{-1}(Y_s)$ in (3.6)–(3.7). The last two observations and (3.6)–(3.7) yield the formulas stated in the lemma. $\square$

Let $\sigma_{ij}(x) = (\partial h_i / \partial x_j) \circ H^{-1}(x)$, where the $h_i$'s and $H$ are as in Lemma 3.1.

LEMMA 3.3. *There exists $c_1$ such that for each $i, j, k \leq d$, and for each $x \in H(B(0, \varepsilon/2) \cap D)$,*

$$\left|\frac{\partial \sigma_{ij}}{\partial x_k}(x)\right| \leq c_1 |x_d|^{\gamma - 1}, \qquad x \in \mathbb{R}^d. \tag{3.8}$$

PROOF. In view of (3.1), all we have to show is that there exists $c_2 > 0$ such that for all $x \in B(0, \varepsilon) \cap D$,

$$h_d(x) \geq c_2 \operatorname{dist}(x, \partial D).$$

But this is immediate from Lemma 3.1(h) and the fact that $h_d = 0$ on $\partial D \cap B(0, \varepsilon)$. $\square$

Let

$$G(x) = |x_d|^{\gamma - 1}. \tag{3.9}$$

PROPOSITION 3.4. *Suppose $\sigma$ maps $\mathbb{R}^d$ to the class of $d \times d$ matrices, $\gamma > 1/2$, and there exists $c_1$ such that for each $i, j, k \leq d$,*

$$\left|\frac{\partial \sigma_{ij}}{\partial x_k}(x)\right| \leq c_1 |x_d|^{\gamma-1}, \qquad x \in \mathbb{R}^d. \tag{3.10}$$

*Consider the stochastic differential equation*

$$dY_t = \sigma(Y_t)\, dW_t + \mathbf{m}(Y_t)\, dL_t^Y, \qquad Y_0 = y_0, \tag{3.11}$$

*where $Y_t$ is a $d$-dimensional process, $W_t$ is a $d$-dimensional Brownian motion with respect to a filtered probability space, $y_0 \in \overline{U}$, $\mathbf{m}$ is the inward pointing normal vector on $\partial U$, and $L_t^Y$ is a local time of $Y$ on $\partial U$. If $Y$ and $Y'$ are two solutions with respect to the same Brownian motion and for some bounded stopping time $\tau$,*

$$\mathbb{E} \int_0^\tau (G(Y_s)^2 + G(Y'_s)^2)\, ds < \infty, \tag{3.12}$$

*then*

$$\mathbb{P}(Y_t = Y'_t \text{ for all } t \in [0, \tau]) = 1.$$

PROOF. Let $\delta > 0$ be small. We can find $c_2, c_3 < \infty$ such that for each $\varepsilon > 0$ there exists a $C^2$ function $f_\varepsilon : \mathbb{R}^d \to \mathbb{R}$ such that $f_\varepsilon(x)$ is a nonincreasing function of $|x|$, $f_\varepsilon(x) = -\log|x|$ for $|x| \geq \varepsilon$, $f_\varepsilon(x) \geq -\log \varepsilon$ if $|x| < \varepsilon$, and

$$|\nabla f_\varepsilon(x)| \leq \frac{c_2}{|x|}, \qquad \left|\frac{\partial^2 f_\varepsilon}{\partial x_i \, \partial x_j}(x)\right| \leq \frac{c_3}{|x|^2}, \qquad i, j \leq d,\ x \in \mathbb{R}^d. \tag{3.13}$$

We see that $\partial f_\varepsilon(x)/\partial x_d \leq 0$ for $x \in \overline{U}$ and $\partial f_\varepsilon(x)/\partial x_d \geq 0$ for $x \in U^c$.

Using (3.10), for any $i, j$ the mean value theorem tells us that

$$\begin{aligned}
|\sigma_{ij}(x) - \sigma_{ij}(y)| &\leq c_4 |x - y|(|x_d|^{\gamma-1} \vee |y_d|^{\gamma-1}) \\
&\leq c_4 |x - y|(G(x) + G(y)).
\end{aligned} \tag{3.14}$$

Let $Z_t = Y_t - Y'_t$ and apply Itô's formula with the function $f_\varepsilon$. We have

$$f_\varepsilon(Z_t) = f_\varepsilon(Z_0) + \int_0^t \sum_{k=1}^d \frac{\partial f_\varepsilon}{\partial x_k}(Z_s)\, dZ_s^k + \frac{1}{2} \int_0^t \sum_{k,\ell=1}^d \frac{\partial^2 f_\varepsilon}{\partial x_k \, \partial x_\ell}(Z_s)\, d\langle Z^k, Z^\ell \rangle_s$$

$$+ \int_0^t \sum_{k=1}^d \frac{\partial f_\varepsilon}{\partial x_k}(Z_s)(\mathbf{1}_{\{k=d\}}\mathbf{m}_k(Y_s) dL_s^Y - \mathbf{1}_{\{k=d\}}\mathbf{m}_k(Y'_s) dL_s^{Y'})$$

$$= f_\varepsilon(Z_0) + \int_0^t \sum_{k=1}^d \frac{\partial f_\varepsilon}{\partial x_k}(Z_s) \left[\sum_{j=1}^d (\sigma_{kj}(Y_s) - \sigma_{kj}(Y'_s))\, dW_s^j\right]$$

$$(3.15) \quad + \frac{1}{2} \int_0^t \sum_{k,\ell=1}^d \frac{\partial^2 f_\varepsilon}{\partial x_k \partial x_\ell}(Z_s)$$

$$\times \left[ \sum_{j=1}^d (\sigma_{kj}(Y_s) - \sigma_{kj}(Y'_s))(\sigma_{\ell j}(Y_s) - \sigma_{\ell j}(Y'_s)) \right] ds$$

$$+ \int_0^t \frac{\partial f_\varepsilon}{\partial x_d}(Z_s)(dL^Y_s - dL^{Y'}_s)$$

$$= f_\varepsilon(Z_0) + M_t + A_t + V_t.$$

Note that if $Y_s \in \partial U$, then $Z_t \in U^c$ and so $\frac{\partial f_\varepsilon}{\partial x_d}(Z_s) \geq 0$, while if $Y'_s \in \partial U$, then $Z_s \in \overline{U}$ and so $\frac{\partial f_\varepsilon}{\partial x_d}(Z_s) \leq 0$. Hence $V_t$ is a nondecreasing process.

By (3.13) and (3.14),

$$|A_\tau| \leq \int_0^\tau \frac{c_5}{|Z_s|^2} |Y_s - Y'_s|^2 (G(Y_s) + G(Y'_s))^2 \, ds$$

$$\leq c_6 \int_0^\tau (G(Y_s)^2 + G(Y'_s)^2) \, ds.$$

Also by (3.13) and (3.14),

$$\langle M \rangle_\tau \leq \int_0^\tau \frac{c_7}{|Z_s|^2} |Y_s - Y'_s|^2 (G(Y_s) + G(Y'_s))^2 \, ds$$

$$\leq c_8 \int_0^\tau (G(Y_s)^2 + G(Y'_s)^2) \, ds.$$

Therefore, using (3.12),

$$\mathbb{E}|A_\tau| \leq c_9, \qquad E \langle M \rangle_\tau \leq c_{10}.$$

So there exists $N > 1$, independent of $\varepsilon$, such that

$$\mathbb{P}(|A_\tau| \geq N) \leq \frac{\mathbb{E}|A_\tau|}{N} \leq \frac{c_9}{N} \leq \delta$$

and

$$\mathbb{P}(|M_\tau| > N) \leq \frac{\mathbb{E}M_\tau^2}{N^2} = \frac{\mathbb{E}\langle M \rangle_\tau}{N^2} \leq \frac{c_{10}}{N^2} \leq \delta.$$

Since $f_\varepsilon(Z_0) \geq \log(1/\varepsilon)$ and $V_t$ is nondecreasing, (3.15) and the above estimates imply that except for an event of probability at most $2\delta$, we have

$$\liminf_{\varepsilon \to 0} f_\varepsilon(Z_\tau) \geq \liminf_{\varepsilon \to 0} f_\varepsilon(Z_0) - 2N = \infty.$$

We conclude that $Z_\tau = 0$, with probability greater than or equal to $1 - 2\delta$. Since $\delta$ is arbitrary, $Z_\tau = 0$, a.s. This is true, a.s., simultaneously for all stopping times $\tau \wedge t$, for rational $t$. Since $Z_t$ is continuous in $t$, our result follows. $\square$

COROLLARY 3.5. *Suppose that $\varepsilon > 0$ is as in Lemma 3.1, $x_0 \in D \cap B(0, \varepsilon/2)$, and $X$ and $X'$ are solutions to*

$$X_t = x_0 + W_t + \int_0^t \mathbf{n}(X_s) \, dL_s^X,$$

$$X'_t = x_0 + W_t + \int_0^t \mathbf{n}(X'_s) \, dL_s^{X'},$$

*relative to the same Brownian motion $W$. Let*

$$\tau = \inf\{t > 0 : X_t \notin D \cap B(0, 3\varepsilon/4) \text{ or } X'_t \notin D \cap B(0, 3\varepsilon/4)\}.$$

*Then*

$$\mathbb{P}(X_t = X'_t \text{ for all } t \in [0, \tau]) = 1.$$

PROOF. Let $H$ be as in Lemma 3.1 and $Y = H(X)$ and $Y' = H(X')$. It will suffice to show that $\mathbb{P}(Y_t = Y'_t \text{ for all } t \in [0, \tau]) = 1$. In view of Lemma 3.2 and Proposition 3.4, all we have to show is that (3.12) holds for these processes $Y$ and $Y'$.

Extend the definition of each $\sigma_{ij}$ by setting $\sigma_{ij}(\hat{x}, x_d) = \sigma_{ij}(\hat{x}, |x_d|)$ for $x \in \mathbb{R}^d \setminus \overline{U}$. By [8] (see also [2]), the Green function $K(x, y)$ for the corresponding elliptic operator is comparable to $|x - y|^{2-d}$. Let $U_1 = H(B(0, \varepsilon) \cap D)$. If $Y_0 = y_0$, then

$$\mathbb{E} \int_0^\tau G(Y_s)^2 \, ds \leq \int_{U_1} K(y_0, x) G(x)^2 \, dx.$$

We will estimate the last quantity in the special case $y_0 = 0$. The general case can be dealt with in an analogous way. Let $A_{kj} = \{x \in U : 2^{-k} \leq |x| < 2^{-k+1}, 2^{-j} \leq x_d < 2^{-j+1}\}$. For $x \in A_{kj}$, we have $K(y_0, x) \leq c_1 2^{-k(2-d)}$ and $G(x)^2 \leq c_2 2^{-2(\gamma-1)j}$. The volume of $A_{kj}$ is bounded by $c_3 2^{-k(d-1)} 2^{-j}$. Since $\gamma > 1/2$, the series $\sum_{j \geq 1} 2^{-j(2\gamma-1)}$ is summable, and we have

$$\int_{U_1} K(y_0, x) G(x)^2 \, dx \leq \sum_{k \geq 1} \sum_{j \geq 1} c_1 2^{-k(2-d)} c_2 2^{-2(\gamma-1)j} c_3 2^{-k(d-1)} 2^{-j}$$

$$\leq c_4 \sum_{k \geq 1} 2^{-k} \sum_{j \geq 1} 2^{-j(2\gamma-1)} < \infty.$$

It follows that $\mathbb{E} \int_0^\tau G(Y_s)^2 \, ds < \infty$ and similarly, $\mathbb{E} \int_0^\tau G(Y'_s)^2 \, ds < \infty$. □

PROOF OF THEOREM 1.1. By standard arguments, the local version of pathwise uniqueness proved in Corollary 3.5 can be extended to the global assertion stated in Theorem 1.1. □

REMARK 3.6. The same proof as that in Proposition 3.4 proves the following theorem.

THEOREM 3.7. *Suppose $Z$ and $Z'$ are two solutions taking values in $\mathbb{R}^d$ to the stochastic differential equation*

$$dZ_t = \sigma_{ij}(Z_t)\,dW_t, \tag{3.16}$$

*where $W_t$ is a $d$-dimensional Brownian motion. (We allow $Z, Z'$ to be adapted to possibly different filtrations.) Suppose there exists a function $G:\mathbb{R}^d \to [0,\infty)$ such that for each $t < \infty$, a.s.,*

$$\int_0^t (G(Z_s)^2 + G(Z'_s)^2)\,ds < \infty,$$

*and for each $i,j = 1,\ldots,d$ and $x,y \in \mathbb{R}^d$, we have*

$$|\sigma_{ij}(x) - \sigma_{ij}(y)| \leq |x-y|(G(x) + G(y)). \tag{3.17}$$

*Then*

$$\mathbb{P}(Z_t = Z'_t \text{ for all } t \geq 0) = 1.$$

REMARK 3.8. The proof in [6] that weak uniqueness holds for (1.1) is correct, but the proof of strong existence is not. If $\{\mathcal{F}_t\}$ is the filtration of the Brownian motion, it was proved that there exists a solution $X_t$ of (1.1) with $X_t$ being $\mathcal{F}_1$ measurable for all $t \leq 1$, but the process $X$ constructed there was not necessarily adapted, that is, it was not shown that $X_t$ was $\mathcal{F}_t$ measurable.

**4. Counterexample.** Let $U$ be the upper half-plane and consider the stochastic differential equation

$$dY_t = \sigma(Y_t)\,dW_t + \mathbf{n}(Y_t)\,dL_t^Y, \qquad Y_0 = y_0, \tag{4.1}$$

where $Y_t$ is a three-dimensional process, $W_t$ is a three-dimensional Brownian motion with respect to a filtered probability space, $y_0 \in \overline{U}$, $\mathbf{n}$ is the inward pointing normal vector and $L_t^Y$ is the local time of $Y$ on $\partial U$.

THEOREM 4.1. *For every $\gamma \in (0, 1/2)$ there exists $\sigma$ which maps $U$ to the class of $3 \times 3$ matrices such that $\sigma$ is bounded, uniformly positive definite, there exists $c_1$ such that for each $i,j,k \leq 3$,*

$$\left|\frac{\partial \sigma_{ij}}{\partial x_k}(x)\right| \leq c_1 x_3^{\gamma-1}, \qquad x \in U, \tag{4.2}$$

*and there exist two solutions $Y$ and $Y'$ to (4.1) with respect to the same Brownian motion and with $y_0 = 0$ such that for some $t > 0$, $\mathbb{P}(Y_t \neq Y'_t) > 0$.*

PROOF. We will identify $\partial U$ and $\mathbb{R}^2$. Choose a function $\Psi : \mathbb{R}^2 \to [0, \infty)$ with the following properties:

(i) $\Psi(x) = \psi(|x|)$ for some $\psi : \mathbb{R} \to [0, \infty)$.
(ii) $\Psi(x)$ is $C^\infty$.
(iii) $\psi(0) = 1$, $\psi$ is decreasing on $[0, \infty)$, $\psi(r) = 0$ for $r > 3/4$.

For integer $k$, let

(4.3) $$\Psi_k^1(x) = 2^{-k\gamma} \Psi(2^k x),$$

(4.4) $$\Psi_k(x) = \sum_{(i_1, i_2) \in \mathbb{Z}^2} \Psi_k^1(x + (3i_1 2^{-k}, 3i_2 2^{-k})).$$

Let $\varphi_k$ be the harmonic function in $U$ which has boundary values $\Psi_k(x)$ on $\partial U$. Consider a large integer $n_1$ whose value will be specified later and let $\varphi = \sum_{m \geq 0} \varphi_{mn_1}$.

Let $\sigma_{ij} \equiv 0$ for $i \neq j$, $\sigma_{33} \equiv 1$ and $\sigma_{jj}(x) = 1 + \varphi(x)$ for $j = 1, 2$.

First we will show that the $\sigma_{ij}$'s satisfy (4.2). For some $c_2, c_3 < \infty$, all $j$, and all $x \in \partial U$, we have $|\varphi_0(x)| \leq c_2$ and $|\frac{\partial \varphi_0}{\partial x_j}(x)| \leq c_3$. Since $\varphi_0$ and $\frac{\partial \varphi_0}{\partial x_j}$ are harmonic, we have $|\varphi_0(x)| \leq c_2$ and $|\frac{\partial \varphi_0}{\partial x_j}(x)| \leq c_3$ for all $x \in U$. By (2.1), $|\frac{\partial \varphi_0}{\partial x_j}(x)| \leq c_4/x_3$ for $x_3 \geq 1$. We use scaling to see that

$$\left| \frac{\partial \varphi_k}{\partial x_j}(x) \right| \leq \begin{cases} c_5 2^{-k(\gamma-1)}, & \text{if } x_3 \leq 2^{-k}, \\ c_6 2^{-k\gamma}/x_3, & \text{if } x_3 > 2^{-k}. \end{cases}$$

Hence, for some $c_7 < \infty$, $j = 1, 2, 3$ and all $x \in U$,

$$\left| \frac{\partial \varphi}{\partial x_j}(x) \right| \leq \sum_{k:x_3 \leq 2^{-k}} c_5 2^{-k(\gamma-1)} + \sum_{k:x_3 > 2^{-k}} c_6 2^{-k\gamma}/x_3 \leq c_7 x_3^{\gamma-1}.$$

This and the definition of the $\sigma_{ij}$'s imply that (4.2) holds.

Assume that pathwise uniqueness holds for (4.1). We will show that this assumption leads to a contradiction. Pathwise uniqueness together with weak existence imply strong existence. Let $Y$ and $Y'$ be the unique solutions to (4.1) with the same driving Brownian motion $W$, with starting points $Y_0$ and $Y_0'$. Most of the time we will be concerned with $Y$ and $Y'$ starting from different points, that is, $Y_0 \neq Y_0'$.

We will write $Y_s = (Y_1(s), Y_2(s), Y_3(s))$, and similarly for $Y_s'$.

Note that the third component of $Y$ is a one-dimensional reflecting Brownian motion starting from 0 and the same is true for $Y'$. Moreover, if $Y_3(0) = Y_3'(0)$, then the third components of $Y$ and $Y'$ are equal to each other for all $t$, a.s.

By Itô's formula, for $k = 1, 2$,

$$(Y_k(t) - Y_k'(t))^2$$
$$= (Y_k(0) - Y_k'(0))^2 + \int_0^t 2(Y_k(s) - Y_k'(s))(\sigma_{kk}(Y_s) - \sigma_{kk}(Y_s')) \, dW_s^k$$

$$+ \tfrac{1}{2} \int_0^t 2(\sigma_{kk}(Y_s) - \sigma_{kk}(Y'_s))^2 \, ds$$

$$= (Y_k(0) - Y'_k(0))^2 + \int_0^t 2(Y_k(s) - Y'_k(s))(\varphi(Y_s) - \varphi(Y'_s)) \, dW^k_s$$

$$+ \int_0^t (\varphi(Y_s) - \varphi(Y'_s))^2 \, ds.$$

It follows that if we set $R_t = |Y_t - Y'_t|^2$, then

(4.5)
$$R_t = R_0 + \int_0^t \sum_{k=1}^2 2(Y_k(s) - Y'_k(s))(\varphi(Y_s) - \varphi(Y'_s)) \, dW^k_s$$
$$+ \int_0^t \sum_{k=1}^2 (\varphi(Y_s) - \varphi(Y'_s))^2 \, ds.$$

This implies that $R_t$ is a time-change of the square of a two-dimensional Bessel process. Hence, the process $t \to |Y_t - Y'_t|$ has the same exit probabilities from an interval as a two-dimensional Bessel process. The main technical goal of this proof will be to show that $R$ hits 1 before some time $t_0 < \infty$ with probability $p_0 > 0$, where $t_0$ and $p_0$ are independent of $R_0 > 0$.

We will now derive some estimates based on (4.5) that will be needed later in the proof.

Suppose that $Y(0) \neq Y'(0)$, $Y_3(0) = Y'_3(0)$ and consider $b > 0$, $a_0 \in [0, 1)$ and $a_1 > 1$. Let

$$\tau_1 = \inf\{t > 0 : |Y_t - Y'_t| \notin [a_0|Y_0 - Y'_0|, a_1|Y_0 - Y'_0|]\},$$
$$\tau_2 = \inf\left\{t > 0 : \int_0^t (\varphi(Y_s) - \varphi(Y'_s))^2 \, ds \geq b|Y_0 - Y'_0|^2\right\}.$$

Since $|Y - Y'|$ is a time-change of a two-dimensional Bessel process, we have for a two-dimensional Brownian motion $\overline{W}$,

$$|Y_t - Y'_t| = \left|Y_0 - Y'_0 + \overline{W}\left(\int_0^t (\varphi(Y_s) - \varphi(Y'_s))^2 \, ds\right)\right|.$$

Note that if $\tau_2 < \infty$, then

$$\int_0^{\tau_2} (\varphi(Y_s) - \varphi(Y'_s))^2 \, ds = b|Y_0 - Y'_0|^2,$$

and therefore

$$|Y_{\tau_2} - Y'_{\tau_2}| = |Y_0 - Y'_0 + \overline{W}(b|Y_0 - Y'_0|^2)|.$$

Hence, if we take $\mathbf{a} = (Y_0 - Y'_0)/|Y_0 - Y'_0|$,

$$\mathbb{P}(\tau_2 < \tau_1) \leq \mathbb{P}(\tau_2 < \infty, |Y_{\tau_2} - Y'_{\tau_2}| \in [a_0|Y_0 - Y'_0|, a_1|Y_0 - Y'_0|])$$
$$\leq \mathbb{P}(|Y_0 - Y'_0 + \overline{W}(b|Y_0 - Y'_0|^2)| \in [a_0|Y_0 - Y'_0|, a_1|Y_0 - Y'_0|])$$
$$= \mathbb{P}(|\mathbf{a} + \overline{W}(b)| \in [a_0, a_1]).$$

The last quantity is a constant $p_0 > 0$ depending on $a_0, a_1$ and $b$ but not depending on $Y$ or $Y'$. Moreover, it is easy to see that for fixed $a_0$ and $a_1$, $p_0 \to 0$ when $b \to \infty$. We record the above inequality for future reference as

$$(4.6) \qquad \mathbb{P}(\tau_2 < \tau_1) \leq p_0.$$

We will now estimate the probability of the complementary event in the same setting. We have

$$\mathbb{P}(\tau_1 < \tau_2) \leq \mathbb{P}\left(\bigcup_{0 \leq t \leq \tau_2} \{|Y_t - Y'_t| \notin [a_0|Y_0 - Y'_0|, a_1|Y_0 - Y'_0|]\}\right)$$

$$\leq \mathbb{P}\left(\bigcup_{0 \leq s \leq b} \{|Y_0 - Y'_0 + \overline{W}(s|Y_0 - Y'_0|^2)| \notin [a_0|Y_0 - Y'_0|, a_1|Y_0 - Y'_0|]\}\right)$$

$$= \mathbb{P}\left(\bigcup_{0 \leq s \leq b} \{|\mathbf{a} + \overline{W}(s)| \notin [a_0, a_1]\}\right).$$

This easily implies that

$$(4.7) \qquad \mathbb{P}(\tau_1 < \tau_2) \leq p_1,$$

where $p_1 \to 0$ as $b \to 0$, for any fixed $a_0$ and $a_1$.

Take $1 - a_0 = a_1 - 1 > 0$ and consider an arbitrary finite stopping time $\tau$. Note that

$$\left\{||Y_\tau - Y'_\tau| - |Y_0 - Y'_0|| > (a_1 - 1)|Y_0 - Y'_0|, \right.$$
$$\left. \int_0^\tau (\varphi(Y_s) - \varphi(Y'_s))^2 \, ds \leq b|Y_0 - Y'_0|^2\right\} \subset \{\tau_1 < \tau_2\}.$$

If we take $a_2 = a_1 - 1$, it follows from (4.7) that for any $a_2, p_0 > 0$ there exists $b > 0$ such that for every stopping time $\tau$,

$$(4.8) \qquad \mathbb{P}\bigg(||Y_\tau - Y'_\tau| - |Y_0 - Y'_0|| > a_2|Y_0 - Y'_0|,$$
$$\int_0^\tau (\varphi(Y_s) - \varphi(Y'_s))^2 \, ds \leq b|Y_0 - Y'_0|^2\bigg) \leq p_0.$$

We will next estimate $\varphi(x) - \varphi(y)$ for $x$ and $y$ such that $x_3 = y_3 > 0$.

The function $\varphi_0$ is not constant on horizontal planes, so for some $x^0 \in U$ we have $\left(\frac{\partial \varphi_0}{\partial x_1}, \frac{\partial \varphi_0}{\partial x_2}\right)(x^0) = \mathbf{v} \neq 0$. Let $r_0 > 0$ be such that $B(x^0, r_0) \subset U$ and

$$(4.9) \qquad \left|\left(\frac{\partial \varphi_0}{\partial x_1}, \frac{\partial \varphi_0}{\partial x_2}\right)(x) - \mathbf{v}\right| \leq |\mathbf{v}|/100$$

for all $x \in B(x^0, r_0)$. By scaling,

$$\left|\left(\frac{\partial \varphi_k}{\partial x_1}, \frac{\partial \varphi_k}{\partial x_2}\right)(x) - 2^{-k(\gamma-1)}\mathbf{v}\right| \leq 2^{-k(\gamma-1)}|\mathbf{v}|/100 \tag{4.10}$$

for $k \geq 1$ and all $x \in B(2^{-k}x^0, 2^{-k}r_0)$.

We have $|\varphi_k(x)| \leq 2^{-k\gamma}$ for all $x$, so $\sum_{k \geq 0} \|\varphi_k\|_\infty < \infty$. It follows that

$$\lim_{n \to \infty} \sup_{x \in B(x^0, r_0)} \sum_{m=1}^{\infty} \varphi_{mn}(x) = 0. \tag{4.11}$$

By (2.1),

$$\lim_{n \to \infty} \sup_{x \in B(x^0, r_0)} \sum_{m=1}^{\infty} |\nabla \varphi_{mn}(x)| = 0.$$

We choose $n_1$ so large that

$$\sup_{x \in B(x^0, r_0)} \sum_{m=1}^{\infty} |\nabla \varphi_{mn_1}(x)| \leq |\mathbf{v}|/100. \tag{4.12}$$

By scaling,

$$\sup_{x \in B(2^{-k}x^0, 2^{-k}r_0)} \sum_{m=k+1}^{\infty} |\nabla \varphi_{mn_1}(x)| \leq 2^{-k(\gamma-1)}|\mathbf{v}|/100. \tag{4.13}$$

For some $c_8 < \infty$, by scaling, $\|\nabla \varphi_k\|_\infty \leq c_8 2^{-k(\gamma-1)}$. We make $n_1$ larger, if necessary, so that

$$\sum_{m=-\infty}^{-1} \|\nabla \varphi_{mn_1}\|_\infty \leq \sum_{m=-\infty}^{-1} c_8 2^{-mn_1(\gamma-1)} \leq |\mathbf{v}|/100. \tag{4.14}$$

We use scaling again to see that

$$\sum_{m=-\infty}^{k-1} \|\nabla \varphi_{mn_1}\|_\infty \leq \sum_{m=-\infty}^{k-1} c_8 2^{-mn_1(\gamma-1)} \leq 2^{-k(\gamma-1)}|\mathbf{v}|/100.$$

Then

$$\sum_{m=0}^{k-1} \|\nabla \varphi_{mn_1}\|_\infty \leq \sum_{m=-\infty}^{k-1} \|\nabla \varphi_{mn_1}\|_\infty \leq 2^{-k(\gamma-1)}|\mathbf{v}|/100. \tag{4.15}$$

With this choice of $n_1$, we have by (4.10), (4.13) and (4.15), for all $x \in B(2^{-k}x^0, 2^{-k}r_0)$,

$$\left|\left(\frac{\partial \varphi}{\partial x_1}, \frac{\partial \varphi}{\partial x_2}\right)(x) - 2^{-k(\gamma-1)}\mathbf{v}\right|$$

$$\leq \left|\left(\frac{\partial \varphi_k}{\partial x_1}, \frac{\partial \varphi_k}{\partial x_2}\right)(x) - 2^{-k(\gamma-1)}\mathbf{v}\right| + \sum_{m=k+1}^{\infty} |\nabla \varphi_{mn_1}(x)| + \sum_{m=0}^{k-1} \|\nabla \varphi_{mn_1}\|_\infty \tag{4.16}$$

$$\leq 3 \cdot 2^{-k(\gamma-1)}|\mathbf{v}|/100.$$

By invariance under rotations by the angle $\pi/2$, there is a vector $\mathbf{w}$, perpendicular to $\mathbf{v}$ and to the vertical axis, such that (4.16) holds with $x^0$ replaced by a different reference point, say, $x^1$; that is, for $k \geq 1$ and $x \in B(2^{-k}x^1, 2^{-k}r_0)$,

$$(4.17) \qquad \left| \left( \frac{\partial \varphi}{\partial x_1}, \frac{\partial \varphi}{\partial x_2} \right)(x) - 2^{-k(\gamma-1)}\mathbf{w} \right| \leq 3 \cdot 2^{-k(\gamma-1)}|\mathbf{w}|/100.$$

Estimates (4.10), (4.13) and (4.15) are invariant under shifts by vectors of the form $(3 \cdot 2^{-k}i_1, 3 \cdot 2^{-k}i_2, 0)$ with $(i_1, i_2) \in \mathbb{Z}^2$, so we have for $(i_1, i_2) \in \mathbb{Z}^2$,

$$(4.18) \qquad \left| \left( \frac{\partial \varphi}{\partial x_1}, \frac{\partial \varphi}{\partial x_2} \right)(x) - 2^{-k(\gamma-1)}\mathbf{v} \right| \leq 3 \cdot 2^{-k(\gamma-1)}|\mathbf{v}|/100,$$

for $x \in B(2^{-k}x^0 + (3 \cdot 2^{-k}i_1, 3 \cdot 2^{-k}i_2, 0), 2^{-k}r_0)$, and

$$(4.19) \qquad \left| \left( \frac{\partial \varphi}{\partial x_1}, \frac{\partial \varphi}{\partial x_2} \right)(x) - 2^{-k(\gamma-1)}\mathbf{w} \right| \leq 3 \cdot 2^{-k(\gamma-1)}|\mathbf{w}|/100,$$

for $x \in B(2^{-k}x^1 + (3 \cdot 2^{-k}i_1, 3 \cdot 2^{-k}i_2, 0), 2^{-k}r_0)$.

Fix some $k_0$ such that $2^{-k_0} < r_0/4$ and consider any vector $\mathbf{u}$ with zero third component such that $2^{-k_0-k-1} \leq |\mathbf{u}| \leq 2^{-k_0-k}$. There exists $c_9 < \infty$ such that for every $z \in \partial U$ we can find $x^2$ and $x^3$ of the form

$$(4.20) \quad \begin{aligned} x^2 &= 2^{-k}x^0 + (3 \cdot 2^{-k}i_1, 3 \cdot 2^{-k}i_2, 0), \\ x^3 &= 2^{-k}x^1 + (3 \cdot 2^{-k}i_1, 3 \cdot 2^{-k}i_2, 0), \end{aligned}$$

with the property that $|z - x^2| \leq c_9 2^{-k}$ and $|z - x^3| \leq c_9 2^{-k}$. One of the vectors $\mathbf{u}$ or $-\mathbf{u}$ must form an angle smaller than $3\pi/8$ with one of the vectors $\mathbf{v}$ or $\mathbf{w}$. Suppose that $\mathbf{u}$ forms an angle smaller than $3\pi/8$ with $\mathbf{v}$ and consider $x \in B(x^2, 2^{-k}r_0/100)$ and $y \in B(x^2 + \mathbf{u}, 2^{-k}r_0/100)$. By integrating $\nabla \varphi$ along the line segment joining $x$ and $y$, and using (4.18), we obtain

$$(4.21) \qquad |\varphi(x) - \varphi(y)| \geq c_{10} 2^{-k\gamma},$$

where $c_{10}$ does not depend on $k, z, x$ or $y$. In the remaining cases, either the above estimate holds as stated or it holds for $x \in B(x^3, 2^{-k}r_0/100)$ and $y \in B(x^3 + \mathbf{u}, 2^{-k}r_0/100)$.

Next we will estimate the rate of growth of $t \to \int_0^t (\varphi(Y_s) - \varphi(Y_s'))^2 \, ds$. Suppose that $Y_0, Y_0' \in \partial U$ and $|Y_0 - Y_0'| \in [2^{-j-1}, 2^{-j}]$. Let $k = j - k_0$, $r_1 = r_0/100$ and $\mathbf{u} = Y_0' - Y_0$. Assume without loss of generality that $\mathbf{u}$ forms an angle smaller than $3\pi/8$ with $\mathbf{v}$. Let $x^2$ be defined as in (4.20), relative to $z = Y_0$. Then (4.21) holds for $x \in B(x^2, 2^{-k}r_1)$ and $y \in B(x^2 + \mathbf{u}, 2^{-k}r_1)$. Let $\tau = \inf\{t > 0 : L_t^Y \geq 2^{-j}\}$. Let $A_1$ be the following event:

The process $Y_3(t)$ reaches level $2^{-j}$ at some time $T_1 < \tau$, then $Y$ hits $B(x^2, r_1 2^{-k-1})$ at a time $T_2 < \tau$, and then stays in $B(x^2, 3r_1 2^{-k}/4)$ until at least $T_2 + 2^{-2k}$.

The probability of $A_1$ is bounded below by $p_1 > 0$, independent of $k$. To see this, we argue as follows. Extend $\sigma$ to all of $\mathbb{R}^3$ by reflection, that is, define $\widetilde{\sigma}(x_1, x_2, x_3) = \widetilde{\sigma}(x_1, x_2, -x_3) = \sigma(x_1, x_2, x_3)$ for $x_3 \geq 0$. Let $\widetilde{V}$ be a weak solution to

$$d\widetilde{V}_t = \widetilde{\sigma}(\widetilde{V}_t) \, dW_t, \qquad \widetilde{V}_0 = (0,0,0). \tag{4.22}$$

The solution to (4.22) is weakly unique since $\widetilde{\sigma}$ is continuous and uniformly positive definite (see [3], Chapter VI). Let $V_t = (V_t^1, V_t^2, V_t^3)$ be defined by $V_t^i = \widetilde{V}_t^i$ for $i = 1, 2$ and $V_t^3 = |\widetilde{V}_t^3|$. A routine calculation with Ito's formula shows that

$$dV_t = \sigma(V_t) \, d\overline{W}_t + \mathbf{n}(V_t) \, dL_t^V$$

for some three-dimensional Brownian motion $\overline{W}$. The proof of the theorem of Yamada–Watanabe (see [10], Section IX.1) applies in the present context, so since we are assuming that there is pathwise uniqueness for (4.1), then there is weak uniqueness for (4.1). Hence $V$ and $Y$ have the same law, and to estimate the probability of $A_1$, it suffices to estimate the probability of $A_1$ when $Y$ is replaced by $V$. By the definition of $V$, it is enough to estimate this probability when $V$ is replaced by $\widetilde{V}$. This latter probability is bounded below, using scaling and the support theorem for diffusions corresponding to uniformly elliptic operators; see [3], Theorem V.2.5. This completes the proof that $\mathbb{P}(A_1) \geq p_1$.

Consider some positive constants $c_{11}$ and $c_{12}$ such that $c_{11} < r_1/8$ and $c_{11} 2^{-k} < |Y_0 - Y_0'|/2$, and let

$$A_2 = \Big\{ ||Y_\tau - Y_\tau'| - |Y_0 - Y_0'|| > c_{11} 2^{-k},$$
$$\int_0^\tau (\varphi(Y_s) - \varphi(Y_s'))^2 \, ds \leq c_{12} 2^{-2k} \Big\}.$$

By (4.8), we can choose $c_{12} > 0$ small so that $\mathbb{P}(A_2) < p_1/4$. Then $\mathbb{P}(A_1 \setminus A_2) > p_1/2$. Suppose that the event $A_1 \setminus A_2$ occurred. Then one possibility is that

$$\int_0^\tau (\varphi(Y_s) - \varphi(Y_s'))^2 \, ds \geq c_{12} 2^{-2k}. \tag{4.23}$$

If the event in (4.23) does not occur, then $Y_t' \in B(x^2 + \mathbf{u}, r_1 2^{-k})$ for $t \in [T_2, T_2 + 2^{-2k}]$ and so by (4.21),

$$\int_0^\tau (\varphi(Y_s) - \varphi(Y_s'))^2 \, ds \geq \int_{T_2}^{T_2 + 2^{-2k}} (\varphi(Y_s) - \varphi(Y_s'))^2 \, ds$$
$$\geq c_{13} 2^{-2k\gamma} 2^{-2k} = c_{13} 2^{-2k(1+\gamma)}. \tag{4.24}$$

Combining the two cases (4.23)–(4.24),

$$\text{(4.25)} \quad \mathbb{P}\bigg(\int_0^\tau (\varphi(Y_s) - \varphi(Y'_s))^2 \, ds \geq c_{14} 2^{-2k(1+\gamma)}\bigg) \geq p_1/2.$$

Recall our assumption that $|Y_0 - Y'_0| \in [2^{-j-1}, 2^{-j}]$. Let

$$T_3 = \inf\{t > 0 : |Y_t - Y'_t| \notin [|Y_0 - Y'_0|/2, 2|Y_0 - Y'_0|]\},$$
$$T_4 = \inf\bigg\{t > 0 : \int_0^t (\varphi(Y_s) - \varphi(Y'_s))^2 \, ds \geq c_{15} 2^{-2j}\bigg\},$$

for some constant $c_{15}$. According to (4.6), we can choose $c_{15}$ sufficiently large so that

$$\text{(4.26)} \quad \mathbb{P}(T_3 > T_4) < 1/2.$$

Recall that $k = j - k_0$ for a suitable constant $k_0$. Let $\tau_0 = 0$, $\tau_1 = \tau \wedge T_3$, $\tau_i = T_3 \wedge (\tau \circ \theta_{\tau_{i-1}} + \tau_{i-1})$ for $i \geq 2$ and

$$F_i = \bigg\{\int_{\tau_i}^{\tau_{i+1}} (\varphi(Y_s) - \varphi(Y'_s))^2 \, ds \geq c_{14} 2^{-2(k+1)(1+\gamma)}\bigg\} \cup \{\tau_{i+1} = T_3\},$$

for $i \geq 0$. The events $F_i$ are not necessarily independent, but if we take $\mathcal{F}_i = \sigma(F_1, \ldots, F_{i-1})$, then we have

$$\mathbb{P}(F_i \mid \mathcal{F}_{i-1}) \geq p_1/2,$$

by (4.25) and the strong Markov property applied at the $\tau_i$'s. Hence, we can compare the sequence $\{F_i\}$ to Bernoulli trials to reach the following conclusion. Let $I_0$ be the smallest integer such that $\sum_{i=0}^{I_0} \mathbf{1}_{F_i} \geq c_{16} 2^{2k\gamma}$ where $c_{16}$ is chosen so that $c_{16} 2^{2k\gamma} c_{13} 2^{-2(k+1)(1+\gamma)} \geq c_{15} 2^{-2j}$. Then $\mathbb{E} I_0 \leq c_{17} 2^{2k\gamma}$. It follows from the definition of $F_i$'s and constant $c_{16}$ that $\tau_{I_0} \geq T_4 \wedge T_3$.

Note that $L^Y_{\tau_{i+1}} - L^Y_{\tau_i} \leq 2^{-j}$. Hence,

$$\text{(4.27)} \quad \mathbb{E} L^Y_{T_4 \wedge T_3} \leq \mathbb{E} L^Y_{\tau_{I_0}} \leq 2^{-j} \mathbb{E} I_0 \leq 2^{-j} c_{17} 2^{2k\gamma} = c_{18} 2^{-k(1-2\gamma)}.$$

Let $S_1 = T_4 \wedge T_3$, $S_i = T_3 \wedge (T_4 \circ \theta_{S_{i-1}} + S_{i-1})$ and $I_1 = \inf\{i : S_i \geq T_3\}$. Let $G_i = \{T_3 \leq S_i\}$ and $\mathcal{G}_i = \sigma(G_1, \ldots, G_i)$. By the strong Markov property applied at the stopping times $S_i$ and (4.26) we have $\mathbb{P}(G_i \mid \mathcal{G}_{i-1}) \geq 1/2$ for all $i$. Thus we can compare the sequence $\{G_i\}$ with a sequence of Bernoulli trials to reach the conclusion that $\mathbb{E} I_1 \leq c_{19} < \infty$.

By (4.27) and the strong Markov property applied at the $S_i$'s, $\mathbb{E}(L^Y_{S_{i+1}} - L^Y_{S_i}) \leq c_{18} 2^{-k(1-2\gamma)}$ for all $i$. A standard argument based on the strong Markov property shows that

$$\text{(4.28)} \quad \begin{aligned} \mathbb{E} L^Y_{T_3} &\leq \mathbb{E} L^Y_{S_{I_1}} \leq c_{18} 2^{-k(1-2\gamma)} \mathbb{E} I_1 \\ &\leq c_{18} 2^{-k(1-2\gamma)} c_{19} = c_{20} 2^{-k(1-2\gamma)}. \end{aligned}$$

Assume that $|Y_0 - Y_0'| = 2^{-j}$ for some (large) integer $j$, and $Y_3(0) = Y_3'(0)$. Let $U_1 = T_3$, $U_i = T_3 \circ \theta_{U_{i-1}} + U_{i-1}$ for $i \geq 2$, and $I_2 = \inf\{i : |Y_{U_i} - Y_{U_i}'| \geq 1\}$. Let $N_k^D$ be the number of downcrossings of the interval $[2^{-k-1}, 2^{-k}]$, that is, $N_k^D$ is the number of distinct $i \leq I_2$ such that $|Y_{U_i} - Y_{U_i}'| = 2^{-k}$ and $|Y_{U_{i+1}} - Y_{U_{i+1}}'| = 2^{-k-1}$. Similarly, let $N_k^U$ be the number of distinct $i \leq I_2$ such that $|Y_{U_i} - Y_{U_i}'| = 2^{-k-1}$ and $|Y_{U_{i+1}} - Y_{U_{i+1}}'| = 2^{-k}$. The process $t \to |Y_t - Y_t'|$ is a time-change of a two-dimensional Bessel process. For any $k \geq 1$, the two-dimensional Bessel process starting at $2^{-k}$ hits 1 before hitting $2^{-k-1}$ with probability $1/(k+1)$, so $N_k^D$ is stochastically majorized by a random variable with a geometric distribution with parameter $1/(k+1)$. Therefore, $\mathbb{E} N_k^D \leq c_{21} k$. Since $|N_k^U - N_k^D| \leq 1$, we also have $\mathbb{E} N_k^U \leq c_{22} k$. Let $N_k$ be the number of distinct $i \leq I_2$ such that $|Y_{U_i} - Y_{U_i}'| = 2^{-k}$. We see that $\mathbb{E} N_k \leq \mathbb{E} N_{k+1}^U + \mathbb{E} N_k^D \leq c_{23} k$.

If $|Y_{U_i} - Y_{U_i}'| = 2^{-k}$, then by (4.28) and the strong Markov property, $\mathbb{E}(L_{U_{i+1}}^Y - L_{U_i}^Y) \leq c_{20} 2^{-k(1-2\gamma)}$. It is routine to combine this with the estimate $\mathbb{E} N_k \leq c_{23} k$ to see that

$$\mathbb{E}\left( \sum_{i \geq 1} (L_{U_{i+1}}^Y - L_{U_i}^Y) \mathbf{1}_{\{|Y_{U_i} - Y_{U_i}'| = 2^{-k}\}} \right) \leq c_{24} k 2^{-k(1-2\gamma)}.$$

Thus $\mathbb{E} L_{U_{I_2}}^Y \leq c_{24} \sum_{k \geq 1} k 2^{-k(1-2\gamma)} < \infty$, due to the assumption $\gamma \in (0, 1/2)$. Let $T_5 = \inf\{t > 0 : |Y_t - Y_t'| = 1\}$ and note that $T_5 = U_{I_2}$. The distribution of $\{L_t^Y, t \geq 0\}$ does not depend on $|Y_0 - Y_0'|$ and we have $L_\infty^Y = \infty$, a.s., so there exist $p_4 > 0$ and $t_0 < \infty$, independent of $|Y_0 - Y_0'|$, such that $\mathbb{P}(U_{I_2} < t_0) = \mathbb{P}(T_5 < t_0) > p_4$. In other words, there exist $p_4 > 0$ and $t_0 < \infty$, independent of $|Y_0 - Y_0'|$, such that $\mathbb{P}(\sup_{t \in [0, t_0]} |Y_t - Y_t'| \geq 1) > p_4$.

We finish by using the following standard argument. We have assumed that there exists a unique strong solution to stochastic differential equation (4.1). Let $Y$ be the solution with $Y_0 = (0, 0, 0)$ and let $Y^j$ be the solution with $Y_0^j = (0, 2^{-j}, 0)$. It is easy to see that the family $\{Y, Y^j, W\}_{\{j \geq 1\}}$ is tight in the space $C([0, \infty), \mathbb{R}^3)^3$ equipped with the topology of uniform convergence on compact sets. Take a weakly convergent subsequence (as $j \to \infty$) and let $\{Y, Y', W\}$ be the weak limit. Note that $(Y^j)^3 = Y^3$ for each $j$, and therefore

$$\int_0^t \mathbf{n}(Y_s^j) \, dL_s^{Y^j} = \int_0^t \mathbf{n}(Y_s) \, dL_s^Y$$

for each $j$ and each $t$. With this observation it is clear that both $Y$ and $Y'$ are solutions to (4.1) with the same driving Brownian motion $W$ and both starting at $(0, 0, 0)$. By the assumption of pathwise uniqueness, $Y_t = Y_t'$ for all $t \geq 0$, a.s. However, we showed that for some $t_0, p_0, \delta > 0$,

$$P\left( \sup_{t \in [0, t_0]} |Y_t - Y_t^j| \geq \delta \right) > p_0$$

for all $j \geq 1$. Passing to the limit, we see that with probability at least $p_0$, we have $\sup_{t \in [0, t_0]} |Y_t - Y'_t| > \delta/2$, a contradiction. Therefore pathwise uniqueness cannot hold for (4.1). $\square$

**Acknowledgments.** We would like to thank G. Lieberman for helpful advice concerning the partial differential equations estimates. We are grateful to the referee for advice on how to simplify some of our arguments and to improve the readability of the paper.

DEPARTMENT OF MATHEMATICS
UNIVERSITY OF CONNECTICUT
STORRS, CONNECTICUT 06269-3009
USA
E-MAIL: bass@math.uconn.edu

DEPARTMENT OF MATHEMATICS
UNIVERSITY OF WASHINGTON
SEATTLE, WASHINGTON 98195
USA
E-MAIL: burdzy@math.washington.edu